\documentclass[11pt]{amsart}
\usepackage[T1]{fontenc}
\usepackage[utf8]{inputenc}
\usepackage{lmodern}
\usepackage[a4paper,left=27mm,right=27mm,top=30mm,bottom=30mm]{geometry}
\usepackage{amsmath,amssymb,amsthm}
\usepackage{graphicx,microtype}
\usepackage[colorlinks=true,linkcolor=blue,citecolor=blue,urlcolor=blue]{hyperref}
\numberwithin{equation}{section}
\newtheorem{theorem}{Theorem}[section]
\newtheorem{lemma}[theorem]{Lemma}
\newtheorem{corollary}[theorem]{Corollary}
\theoremstyle{remark}

\newcommand{\R}{\mathbb R}
\newcommand{\Sn}{\mathbb S^{n-1}}
\newcommand{\ip}[2]{\langle #1,#2\rangle}

\title[Convergence rate of RMCF over asymptotically conical shrinkers]{Convergence rate of rescaled mean curvature flows over asymptotically conical shrinkers}
\author{Qiyuan Song}
\begin{document}

\begin{abstract}
We study how rapidly a rescaled mean curvature flow can approach an asymptotically conical self-shrinker. For a global normal graph with small $C^1$ norm and time-integrable $C^0$ norm, we prove that every nonzero solution has a leading asymptotic term given by an eigenfunction of the linearized operator $L$. 
In particular, a global graph $u$ with $e^{kt}\|u\|_{C^1}\to 0$ for every $k$ must coincide with the shrinker.
We also construct a nonzero local graph over an immersed asymptotically conical shrinker with $C^2$ super-exponential decay in every bounded subset. This flow is defined on expanding open graphical domains. The example illustrates how rigidity fails in the local graphical setting. 
\end{abstract}
\maketitle

\section{Introduction}

Mean curvature flow is a family of immersions $X: M^n\times I\to \mathbb R^{n+1}$ evolving by the mean curvature vector:
\[
 \partial_t X=\vec{H},
\]
where $\vec{H}$ is the mean curvature vector of the hypersurface $M_t=X(M,t)$.

If a mean curvature flow develops a singularity at $(0, 0)\in\mathbb R^{n+1}\times \mathbb R$, we can do a rescaling around it to get the rescaled flow $\tilde M_t=e^{t/2}M_{-e^{-t}}$, which satisfies the equation
\begin{equation}\label{eq:rmcf}
    \ip{\partial_t \tilde X}{\nu}=-H+\frac{1}{2}\ip{\tilde X}{\nu},
\end{equation}
where $\nu$ is the normal vector, and $\vec{H}=-H\nu$.

A self-shrinker is a hypersurface $\Sigma$ satisfying $H_\Sigma=\frac12\ip{X}{\nu_\Sigma}$, which is a stationary solution of rescaled mean curvature flow. 
%and generates the homothetically shrinking flow $M_t=\sqrt{-t}\,\Sigma$,  $t<0$. 
Self-shrinkers provide models for singularities of mean curvature flow.
Once a rescaled flow converges to such a model, it is natural to ask how
rapidly this convergence can occur. In particular, can a nonstationary
flow converge to the shrinker faster than every exponential of the rescaled time?

For compact singularity models, Martin-Hagemayer and Sesum
\cite{martinhagemayer2025rate} proved that a rescaled mean curvature flow of closed hypersurfaces converging to a compact self-shrinker faster than every exponential must coincide with that shrinker.  
Later, Daniels-Holgate and Hershkovits
\cite{danielsholgate2025closeclosesingularmean} extended the result to comparison between two flows. They proved that two compact mean curvature flows with the same multiplicity-one compact singularity must agree if their distance decays faster than every power of time. After rescaling, this becomes a super-exponential closeness condition. 

For noncompact shrinkers,  Huang and Zhao \cite{huangzhao2025rate} proved global graphs over cylinders with small $C^1$ norm and super-exponential $L^2$ decay must be identically $0$. However, the rigidity fails without the global graph assumption. They constructed a Tikhonov-type local graphical solutions with super-exponential  decay on expanding domains.

In this paper, we study the analogous question for asymptotically conical shrinkers. A shrinker $\Sigma\subset\mathbb R^{n+1}$ is called asymptotically conical if there is a regular cone $C$ with vertex at origin, such that
\[
 \lambda\Sigma\longrightarrow C
 \quad\text{in }C^\infty_{\mathrm{loc}}(\R^{n+1}\setminus\{0\})
 \quad\text{as }\lambda\to0^+
\]
Uniqueness of asymptotically conical tangent flows was established by Chodosh and Schulze \cite{chodoshschulze2019uniqueness}.
Chen, Lee, Sun, and Zhao \cite{chen2026closedmeancurvatureflows} showed that every embedded asymptotically conical shrinker can be realized as a multiplicity-one tangent flow of a closed embedded mean curvature flow. 
They also constructed such flows with prescribed asymptotics of the rescaled flow, given by eigenfunctions of the linearized operator $L$.
Our first result describes the asymptotic behaviour under the global graphical assumption. 
\begin{theorem}\label{main thm}
Let $\{\tilde M_t\}_{t\in[0, \infty)}$ be a smooth rescaled mean curvature flow, which is a graph $u$ over an embedded asymptotically conical shrinker $\Sigma$. For sufficient small $\varepsilon>0$, if $u$ satisfies
\begin{equation}\label{cdt1}
\sup_{t\geq0}\|u(t)\|_{C^1(\Sigma)}<\varepsilon,\qquad \int_0^\infty\|u(t)\|_{C^0(\Sigma)}dt<\infty.
\end{equation}
Then either $u\equiv0$, or there are an eigenvalue $\lambda_*$ of $L$
and a nonzero eigenfunction $\varphi_*$ such that
 \begin{equation}\label{case 2}
 e^{-\lambda_*t}u(\cdot,t)\longrightarrow\varphi_*
 \quad\text{in }L^2(\Sigma,d\tilde\mu),\qquad
 L\varphi_*=\lambda_*\varphi_*.
\end{equation}
The convergence also holds in $W^{m,2}(\Sigma,d\tilde\mu)$ for every $m\in\mathbb N$ and $C^\infty_{loc}(\Sigma)$.
\end{theorem}

Here the weighted measure $d\tilde \mu$ and linearized operator $L$ is defined by
\begin{equation}\label{eq:gaussian-L}
 d\tilde\mu=e^{-|X|^2/4}\,d\mu,\qquad
 Lf=\Delta_\Sigma f-\frac12\ip{X}{\nabla_\Sigma f}
       +(|A_\Sigma|^2+\frac12)f
\end{equation}
Throughout the paper, $\|\cdot\|_{L^2}$ and $\ip{\cdot}{\cdot}$ refers to the $L^2$ norm in this weighted space and its inner product. 
All $C^k$ norms are intrinsic norms on the base, and
$\|f\|_{W^{m,2}}^2=\sum_{j=0}^m\|\nabla^j f\|_{L^2}^2$.

As a corollary, we get the rigidity under super-exponential decay.
\begin{corollary}\label{cor:rigidity}
Let $\{\tilde M_t\}_{t\in[0, \infty]}$ be a smooth rescaled mean curvature flow, which is a graph $u$ over an embedded asymptotically conical shrinker $\Sigma$. 
For sufficient small $\varepsilon>0$, if $u$ satisfies
\begin{equation}\label{cdt2}
       \sup_{t\geq 0}\|u(\cdot, t)\|_{C^1(\Sigma)}\leq\varepsilon,\qquad
       \lim_{t\to \infty} e^{kt}\|u(\cdot, t)\|_{C^0(\Sigma)}=0 \text{ for every $k\in\mathbb{N}$}.
    \end{equation} 
Then $u\equiv0$.
\end{corollary}

Our second result constructs a nonzero local graph over the immersed shrinker constructed by Drugan and Kleene \cite{drugan2013immersedselfshrinkers}.
The flow is defined on expanding open domains and has $C^2$ super-exponential convergence on every bounded domain. This example shows that the global super-exponential convergence assumption cannot be weakened to local, at least in the immersed setting.
\begin{theorem}\label{counterexample}
For each $n\geq2$, there is an immersed, rotationally
symmetric asymptotically conical shrinker
\[
 X:\Sigma=\R\times\Sn\longrightarrow\R^{n+1},\qquad
 X(s,\omega)=(z(s),r(s)\omega),
\]
where $s$ is arclength and $r(s)>0$, with the following property.
There exist constants $a,C,T_0>0$ and a smooth rotationally symmetric function $u(s,t)$ on $t>T_0$, $|s|<at$, such that $X+u\nu$ satisfies the rescaled mean curvature flow equation \eqref{eq:rmcf}, with initial data
\begin{equation}\label{initial data}
 u(0,t)=e^{-t^2},\qquad u_s(0,t)=0,
\end{equation}
and
\begin{equation}\label{exp decay}
 \sup_{|s|<at}
 \bigl(|u|+|\nabla_\Sigma u|+|\nabla_\Sigma^2u|\bigr)
 \leq C e^{-t^2/4}.
\end{equation}
Consequently, for every bounded parameter region $V\subset\Sigma$
and every $k>0$,
\[
 \lim_{t\to\infty}e^{kt}\|u(\cdot,t)\|_{C^2(V)}=0.
\]
Every such region is in the graph domain for all sufficiently large $t$. The same applies to $X^{-1}(B_R)$ for each fixed $R>0$.
\end{theorem}

For the first result, we first use the frequency argument in Daniels-Holgate and Hershkovits's work \cite{danielsholgate2025closeclosesingularmean} to give an upper bound for $V=\ip{u}{Lu}/\|u\|^2_{L^2}$. We also show that $V$ converges to an eigenvalue $\lambda_*$ of $L$. The discrete spectrum of $L$ and time-integrable error forces the projection to other eigenspaces goes to $0$ in $L^2$, and $e^{-\lambda_*t}u$ converges to an eigenfunction of $\lambda_*$. The $L^2$ convergence can be improved to $W^{m, 2}$ by energy estimates and compactness of weighted spaces. The time-integrable $C^0$ norm condition can be improved to $C^2$, which is used to control the error term in the rescaled equation.

For the second result, we rewrite the graph equation as a first-order ODE system, and expand in time near a large center $T$. We first solve the finite truncated equations, and give a uniform estimate on the coefficient norm. Passing to subsequence, these solutions will converge to the actual solution. The solutions centered in different time agree on overlaps, so they can be joined into one flow. The super-exponential decay comes from the initial data and a suitable choice of domain.

The paper is organized as follows: Section~\ref{sec:global} proves the global theorem \ref{main thm}, and
Section~\ref{sec:local} constructs the local example \ref{counterexample}.
The calculations of the profile curve in the example are collected in Appendix~\ref{sec:appendix}.

The author used ChatGPT to assist in the construction of the local example. The author checked all mathematical statements and proofs, and takes full responsibility for the content.

\section{Asymptotic behavior of global graphs}\label{sec:global}

As in \cite{danielsholgate2025closeclosesingularmean}, we define the following quantities:
\begin{equation}\label{eq:global-quantities}
\begin{aligned}
 E&=(\partial_t-L)u,\qquad I(t)=\|u(\cdot,t)\|_{L^2}^2,\\
 V(t)&=\frac{\ip{u}{Lu}}{I(t)}
 =\frac{-\|\nabla u\|_{L^2}^2
       +\int_\Sigma(|A|^2+\tfrac12)u^2\,d\tilde\mu}{I(t)}.
\end{aligned}
    \end{equation}

Throughout this section, $\Sigma, u, \varepsilon$ are as in Theorem ~\ref{main thm}.

\begin{lemma}\label{reg lem}
Let $u$ be a smooth global graphical rescaled flow over $\Sigma$ satisfying \eqref{cdt1}.
Then for every $k\in \mathbb N$, there exists constant $C(\Sigma, k)$ such that
\begin{equation}\label{eq:delayed-C2}
 \sup_{t\in[s+\tfrac 1 2,s+1]}\|u(t)\|_{C^k(\Sigma)}\leq C(k,\Sigma) \|u(s)\|_{C^0(\Sigma)},\qquad s\geq 1
\end{equation}
Moreover, $\|u(t)\|_{C^k(\Sigma)}\in L^1\cap L^\infty([2,\infty))$.
\end{lemma}
\begin{proof} 

By \cite[Proposition 4.1]{chodoshchoischulze2023genericII}, sufficient small $C^1$ norm of $u$ and asymptotically conical structure give
\begin{equation}\label{eq:bootstrap}
    \sup_{t\geq 1}\|u(t)\|_{C^m(\Sigma)}\leq C_m<\infty,\qquad \text{for every }m\in\mathbb N.
\end{equation}

The normal-graph equation for the rescaled flow is $\partial_tu=Lu+E(x, u,\nabla u, \nabla^2u)$.
The error term $E(x,z,p,Q)$ is smooth and satisfies $E(x,0,0,Q)=0$. Integrating its first derivatives in $(z,p)$ with fixed $Q=\nabla^2u$, we get $E(x,u,\nabla u,\nabla^2 u)=a_u u+\langle b_u,\nabla u\rangle$. 
The normal-graph formula \cite[Appendix A]{chodoshchoischulze2023genericII} and \eqref{eq:bootstrap} guarantee that $a_u,b_u$ and their spatial derivatives are bounded on $t\geq1$. Hence we can rewrite the equation as
\begin{equation}\label{eq:fixed-graph}
 \partial_tu=\Delta_\Sigma u-\langle{\tfrac{1}{2}X^T-b_u}, \nabla u\rangle+(|A_\Sigma|^2+\frac{1}{2}+a_u)u
\end{equation}

Commuting covariant derivatives in \eqref{eq:fixed-graph}, we get
\begin{equation}\label{eq uj}
    \partial_t \nabla^j u=(\Delta_\Sigma-\nabla_{X^{T}/2-b_ u})\nabla^j u+(|A_\Sigma|^2+\frac{1}{2}+a_u)\nabla^j u+R_j
\end{equation}
where
 \[
 R_j=\sum_{l=0}^{j-1} \nabla^l \mathrm{Rm}*\nabla^{j-l} u\
 +\sum_{l=1}^{j}\nabla^l X*\nabla^{j+1-l} u
    +\sum_{l=0}^{j-2} X*\nabla^l\mathrm{Rm}*\nabla^{j-1-l}u
 \]
 The curvature estimate for asymptotically conical shrinkers \cite[Lemma 2.1]{wang2011uniquenessselfsimilarshrinkersasymptotically} gives
\begin{equation}\label{nabla A}
    |\nabla^j A_\Sigma|\leq C_j(1+|X|)^{-1-j},
 \quad j\geq0
\end{equation}
Since $\mathrm{Rm}_{ijkl}=A_{ik}A_{jl}-A_{il}A_{ij}$, we see that $|\nabla^l X|$, $(1+|X|^{l+2})|\nabla^l \mathrm{Rm}|$ are all bounded. Hence $|R_j|\leq C\sum_{0\leq l\leq j}|\nabla^l u|$ for some constnt $C$.
 
 Set $D=\partial_t-\Delta_\Sigma+\nabla_{X^T/2-b_u}$. Multiplying 
 $\nabla^j u$ to both sides of \eqref{eq uj}, integration by parts and Young's inequality implies
\[
 D|\nabla^j u|^2
 \leq -2|\nabla^{j+1}u|^2+C\sum_{l=0}^j|\nabla^l u|^2
\]

Fix $s>1$, for $s\leq t\leq s+1$, let
\[
 F(x,t)=\sum_{j=0}^k\frac{(t-s)^j}{2^j j!}|\nabla^j u|^2
\]
The negative gradient terms will absorb the derivatives of the
time factors, hence $DF\leq C F$ for some constant $C$.
The shrinker identity gives
$D(1+|X|^2)=-2n+|X|^2-2\langle b_u,X^T\rangle\geq-K(1+|X|^2)$, for some constant $K$. For $\delta>0$, let
\[
G(x,t)=e^{-C(t-s)}F-\delta e^{K(t-s)}(1+|X|^2),
\]
then $DG\leq 0$. We apply maximum principle to $G$ in $X^{-1}(B_R(0))\times [s, s+1]$.

Along the parabolic boundary, $G(\cdot,s)\leq \|u(s)\|^2_{C^0}$, $G|_{\{|X|=R\}}\leq \sup_{t\geq 1}\|u(t)\|_{C^k}^2-\delta R^2<0$ for $R>C_k^2/\delta$. By maximum principle, $e^{-C(t-s)}F(x,t)\leq\|u(s)\|^2_{C^0}+\delta e^{K(t-s)}(|X|^2+1)$ in $X^{-1}(B_R(0))$.
Fix $x$, sending $\delta\to 0^+$
gives $F(x,t)\leq e^{C(t-s)}\|u(s)\|_{C^0}^2$.
For $s+1/2\leq t\leq s+1$, we have $\sum_{0\leq j\leq k} |\nabla^j u|^2(x,t)\leq4^{k}k!F(x,t)$. This proves \eqref{eq:delayed-C2}.
Integrating this on $t-1\leq s\leq t-\tfrac{1}{2}$ yields
\[
\|u(t)\|_{C^k}\leq 2C\int_{t-1}^{t-\tfrac{1}{2}}\|u(s)\|_{C^0}ds
\]
By Tonelli's theorem,
\[
\int_2^\infty\int_{t-1}^{t-\tfrac{1}{2}}\|u(s)\|_{C^0}dsdt
\leq2\int_1^\infty\int_{s+\tfrac{1}{2}}^{s+1}\|u(s)\|_{C^0}dtds=\int_1^\infty \|u(s)\|_{C^0}ds<\infty
\]
Hence $\|u(t)\|_{C^k}\in L^1\cap L^\infty([2,+\infty))$.

\end{proof}

For the remainder of the section, set $c(t)=\|u(t)\|_{C^2}$. Shifting the initial time if necessary, we can assume $c\in L^1\cap L^\infty([0,\infty))$.

\begin{lemma}\label{E le}
There exists some constants $\Lambda, C$ such that $\sup_{t\geq 0}V(t)\leq \Lambda$, 
\begin{equation}\label{E1}
 \|E\|_{L^2}^2\leq Cc(t)^2I(t)(\Lambda-V(t)).
\end{equation}
\end{lemma}
\begin{proof}
Let $a=\sup_\Sigma |A|^2+\frac{1}{2}$, then $V\leq a$,
\begin{equation}\label{eq Du}
    \int_\Sigma |\nabla u|^2d\tilde \mu=\int_\Sigma (|A|^2+\frac{1}{2})u^2d\tilde\mu-VI    \leq (a-V)I
\end{equation}

By the normal-graph expansion 
\cite[Appendix A]{danielsholgate2025closeclosesingularmean} and the curvature estimates\eqref{nabla A}, there exists a constant $C$ depending on $\Sigma, \varepsilon$ such that
\[
 |E|\leq C\bigl(|u|+|\nabla u|+|\nabla^2u|\bigr)
                 \bigl(|u|+|\nabla u|\bigr).
\]
Since $\|u(t)\|_{C^2}=c(t)\in L^1\cap L^\infty([0,\infty))$, we have
\begin{align*}
    \|E\|^2_{L^2}&\leq C^2c^2(t)\int_\Sigma u^2+|\nabla u|^2d\tilde \mu\\
    &\leq C^2c^2(t)(I+(a-V)I)\leq C^2c^2(t)I(\Lambda-V)
\end{align*}
Here $\Lambda=1+a>\sup_{t\geq 0} V(t)$.

\end{proof}

If $u\equiv0$, then there is nothing to prove. Now we assume $u$ is not identically $0$ and prove \eqref{case 2}. The next two lemmas show that $V$ is bounded from below and has a limit.

\begin{lemma}\label{V lem}
There exists $V_0>-\infty$, such that $V(t)\ge V_0$ for all $t$.
\end{lemma}
\begin{proof}
We compute the derivative of $I$ and $V$. Self-adjointness of $L$ gives
$I'=2\ip{u}{Lu+E}$,
\begin{align}
        V'&=\frac{1}{I}(\langle \partial_tu, Lu\rangle+\langle u, L\partial_tu\rangle)-V\frac{I'}{I} \notag\\
        &=\frac{2}{I}\langle \partial_tu, Lu\rangle-2V(V+\frac{\langle u, E\rangle}{I})\notag\\
        %&=\frac{2}{I}(\langle Lu, Lu\rangle+\langle Lu, E\rangle-\langle Vu, Vu\rangle-\langle Vu, E\rangle)\notag\\
        &=\frac{2}{I}(\|(L-V)u\|^2_{L^2}+\langle E,(L-V)u\rangle)
        \geq -\frac{\|E\|^2}{2I} \label{eq V'}
    \end{align}
Hence $(\Lambda-V)'\leq Cc(t)^2(\Lambda-V)$. Gronwall inequality implies
\[
 V(t)\geq\Lambda-(\Lambda-V(0))
       \exp\!\left(C\int_0^t c^2(\tau)\,d\tau\right).
\]
$c^2$ is integrable as $c\in L^1\cap L^\infty([0,\infty))$, hence $V$ is bounded from below for all time. 
\end{proof}

\begin{lemma}\label{lambda lem}
There exists an eigenvalue $\lambda_*$ of $L$ such that $V(t)\to\lambda_*$
as $t\to\infty$.
\end{lemma}
\begin{proof}
Completing the square in \eqref{eq V'} gives
\begin{align*}
       V'&\geq \frac{2}{I}(\|(L-V)u\|^2_{L^2}+\langle E,(L-V)u\rangle)\\
       &\geq \frac{\|(L-V)u\|^2_{L^2}}{I}-Cc^2(t)(\Lambda-V)
   \end{align*}
Let $\tilde V(t)=V(t)+\int_0^t Cc^2(s)(\Lambda-V(s))ds$, then
   \[
   \tilde V'=V'+Cc^2(t)(\Lambda-V)\geq \frac{\|(L-V)u\|^2_{L^2}}{I}\geq 0
   \]
   %\[
   %\tilde V(t)=V(t)+\int_0^t Cc^2(s)(\Lambda-V(s))ds\leq \Lambda %+\int_0^t C\Lambda c^2(s)ds<\Lambda(1+C\|c\|^2_{L^2})
   %\]
$\tilde V$ is nondecreasing and bounded from above, hence has a limit. Since $c^2$ is integrable, $V$ also has a limit $\lambda_*$.
Note that $S=\|(L-V)u\|^2_{L^2}/I$ is integrable in time:
\[
\int_0^\infty \frac{\|(L-V)u\|^2_{L^2}}{I}dt \leq \int_0^\infty \tilde V'(t)dt=\tilde V_*-\tilde V(0)< \infty
\]
Hence we can take a sequence $t_i\to\infty$, such that $S(t_i)\to 0$. Let $w=u/\|u\|_{L^2}$, then $S(t_i)=\|(L-V)w(t_i)\|^2_{L^2}\to 0$.
By spectrum theory, $\text{dist}^2(V(t_i), \sigma(L))\leq \|(L-V)w(t_i)\|^2_{L^2}$. Sending $i\to \infty$, we get $\text{dist}(\lambda_*, \sigma(L))=0$. Since $\sigma(L)$ is discrete, this implies $\lambda_*\in \sigma(L)$.
   
\end{proof}

\begin{lemma}\label{lem:projection}
Let $P_*$ be the orthogonal projection onto $E_{\lambda_*}=\ker(L-\lambda_*)$.
For $w=u/\|u\|_{L^2}$,
\begin{equation}\label{w-Pw}
 \lim_{t\to\infty}\|w-P_*w\|_{L^2}=0
\end{equation}
\end{lemma}
\begin{proof}
The evolution equation for $w$ is
\begin{equation}\label{eq w}
    \partial_t w=(L-V)w+\frac{E-w\langle E, w\rangle}{\|u\|_{L^2}}=(L-V)w+F
\end{equation}
 Choosing larger $C$ if necessary, equation \eqref{E1} and Lemma ~\ref{V lem} implies $\|E\|_{L^2}\leq Cc(t)\|u\|_{L^2}$. Then
\[
\|F\|_{L^2}=\frac{\|E-w\langle E, w\rangle\|}{\|u\|_{L^2}}
\leq \frac{2\|E\|_{L^2}}{\|u\|_{L^2}}\leq 2Cc(t)
\]

By \cite[Lemma 3.1]{chodosh2024meancurvatureflowgeneric}, $L$ has eigenvalues $\lambda_1> \lambda_2\geq\lambda_3\geq\dots $, with $\lim_{i\to\infty}\lambda_i=-\infty$, and there exists an orthonormal basis $\{\varphi_i\}$ of $L^2(\Sigma)$, such that $L\varphi_i=\lambda_i\varphi_i$.
Let $E_-=\text{span}\{\varphi_i: \lambda_i<\lambda_*\}$, $E_+=\text{span}\{\varphi_i: \lambda_i>\lambda_*\}$, and let $P_-, P_+$ be the projections to the corresponding eigenspaces. We will prove $\|P_-w\|_{L^2}, \|P_+w\|_{L^2}$ both converge to $0$.

Since $\sigma(L)$ is discrete,there exists $b>0$ such that
$\inf_{\lambda_i\neq \lambda_*}|\lambda_i-\lambda_*|> 2b$.
By Lemma ~\ref{lambda lem}, there exists $T>0$ such that $|\lambda_i-V(t)|\geq b$ for $\lambda_i\neq \lambda_*$, $t\geq T$. This implies $\ip{(L-V)P_-w}{P_-w}\leq -b\|P_-w\|_{L^2}^2$, $\ip{(L-V)P_+w}{P_+w}\geq b\|P_+w\|_{L^2}^2$, for $t\geq T$.

Projecting the equation\eqref{eq w} to $E_\pm$ and taking inner products with $P_\pm w$, we get
\[
 \frac d{dt}\|P_-w\|_{L^2}\leq-b\|P_-w\|_{L^2}+Cc(t),\qquad
 \frac d{dt}\|P_+w\|_{L^2}\geq b\|P_+w\|_{L^2}-Cc(t).
\]
For $t>T$, 
\[
\|P_-w(t)\|_{L^2}\leq e^{-b(t-T)}\|P_-w(T)\|_{L^2}+C\int_T^t e^{-b(t-\tau)}c(\tau)d\tau
\]
The first term goes to $0$ as $t\to\infty$. For any $\delta>0$, there exists $T'>T$ such that $\int_{T'}^\infty c(s)ds\leq \delta$. 
Then for $t\geq T'$ large enough such that $e^{-b(t-T')}\leq \delta/\|c\|_{L^1}$, we have
\begin{equation*}
    \int_T^te^{-b(t-\tau)}c(\tau)d\tau\leq 
    \int_T^{T'}e^{-b(t-T')}c(\tau)d\tau+\int_{T'}^\infty c(\tau)d\tau\leq 2\delta
\end{equation*}
Sending $t\to\infty$, $\delta\to 0$, we get $\|P_-w(t)\|_{L^2}\to0$.

For $s>t>T$,
\[
 \|P_+w(t)\|_{L^2}
 \leq e^{-b(s-t)}\|P_+w(s)\|_{L^2}+C\int_t^s e^{-b(\tau-t)}c(\tau)\,d\tau.
\]
Sending $s\to\infty$, we get
\[
\|P_+w(t)\|_{L^2}\leq 0+C\int_t^\infty e^{-b(\tau-t)}c(\tau)d\tau\leq C\int_t^\infty c(\tau)d\tau
\]
Thus $\|P_+w(t)\|_{L^2}\to 0$ when $t\to\infty$. Then $\|w-P_*w\|_{L^2}=\|P_-w+P_+w\|_{L^2}\to 0$.

\end{proof}

\begin{proof}[Proof of Theorem \ref{main thm}]
    Assume $u$ is not identically $0$. We first show the $L^2$ convergence.
    Let $z=e^{-\lambda_* t}P_* u$, then
    \[
    z'=e^{-\lambda_* t}(P_*(Lu+E)-\lambda_*P_*u)=e^{-\lambda_* t}P_*E
    \]
    By \eqref{w-Pw}, we can choose large $T$ such that $ \|u\|_{L^2}\leq 2\|P_*u\|_{L^2}$ for $t\geq T$. Then
    \begin{align*}
       | \tfrac{d}{dt}\|z\|_{L^2}|
        %=\frac{\langle z', z\rangle}{\|z\|_{L^2}}
        &\leq \| z'\|_{L^2}=\|e^{-\lambda_* t}P_*E\|_{L^2}\\
        &\leq e^{-\lambda_* t}\|E\|_{L^2}\leq Cc(t)e^{-\lambda_* t}\|u\|_{L^2}\leq 2Cc(t)\|z\|_{L^2}
    \end{align*}
    Gronwall inequality gives
   \[
   \|z(T)\|_{L^2}\exp{(-\int_T^t Cc(t)dt)}\leq \|z(t)\|_{L^2}\leq \|z(T)\|_{L^2}\exp{(\int_T^t Cc(t)dt)}
   \]
   Since $u$ is not identically $0$, we have $\inf_{t\geq T} \|z(t)\|_{L^2}>0$.  For any $t>s>T$, 
    \[
    \|z(t)-z(s)\|_{L^2}\leq \int_s^t \|z'(\tau)\|_{L^2}d\tau
    \leq \int_s^t Cc(\tau)\|z(\tau)\|_{L^2}d\tau
    \leq C\sup_{t\geq T} \|z(t)\|_{L^2}\int_s^t c(\tau)d\tau
    \]
    The integrability of $c(t)$ then implies $z(t)$ is Cauchy in $E_{\lambda_*}$, so $z(t)$ has a $L^2$ limit $\varphi_*\in E_{\lambda_*}$. $\varphi_*$ is nonzero since $\inf_{t\geq T} \|z(t)\|_{L^2}>0$.
    
    By \eqref{w-Pw},  $\ \|u-P_*u\|_{L^2}\leq \frac{1}{2}\|u\|_{L^2}$ for sufficiently large $t$. Thus
    \[
    \lim_{t\to\infty} \|e^{-\lambda_* t} (u-P_* u)\|_{L^2}
    =\lim_{t\to\infty} \frac{\|u-P_*u\|_{L^2}}{\|P_*u\|_{L^2}}\|e^{-\lambda_* t}P_*u\|_{L^2}
    \leq\lim_{t\to\infty} \frac{\|u-P_*u\|_{L^2}}{\frac{1}{2}\|u\|_{L^2}}\|z\|_{L^2}=0
    \]
    This proves the $L^2$ convergence. Next, we improve it to $W^{m, 2}$ convergence.

Set $\tilde u=e^{-\lambda_*t}u$, the evolution equation for $\tilde u$ is
\[
\partial_t\tilde u=(\Delta_\Sigma-\nabla_{X^{T}/2-b_u})\tilde u+(|A_\Sigma|^2+\frac{1}{2}-\lambda_*+a_u)\tilde u\]
Similar to the process in Lemma~\ref{reg lem}, we can get
\begin{equation}\label{eq Cj}
 \frac d{dt}\|\nabla^j\tilde u\|_{L^2}^2+\|\nabla^{j+1}\tilde u\|_{L^2}^2\leq C_j\sum_{k=0}^j\|\nabla^k \tilde u\|_{L^2}^2.
\end{equation}
 Fix $s>1$, let
\[
U_m(t)=\sum_{j=0}^m \frac{(t-s)^j}{2^jj!}\|\nabla^j\tilde u\|_{L^2}^2
\]
For $t\in[s, s+1]$,
\begin{align*}
\frac{d}{dt} U_m&\leq -\sum_{j=0}^m \frac{(t-s)^j}{2^{j+1}j!}\|\nabla^{j+1}\tilde u\|_{L^2}^2
+\sum_{j=0}^{m}C_j\frac{(t-s)^j}{2^jj!}\sum_{k=0}^j\|\nabla^{k}\tilde u\|_{L^2}^2\\
&\leq 0+\hat{C}_m\sum_{k=0}^{m}\big(\sum_{j=k}^m\frac{(t-s)^j}{2^j}\big)\|\nabla^{k}\tilde u\|_{L^2}^2\leq 2\hat{C}_mU_m
\qquad \hat{C}_m=\max_{1\leq j\leq m}\{C_j\}
\end{align*}
Since $U_m(s)=\|\tilde u(\cdot, s)\|_{L^2}^2$, $U_m(s+1)\geq (2^m m!)^{-1}\|\tilde u(\cdot, s+1)\|_{W^{m, 2}}^2$, we have
\begin{equation}\label{Um}
    \|\tilde u(\cdot, s+1)\|_{W^{m, 2}}\leq  e^{2\hat{C}_m}2^m m!\|\tilde u(\cdot, s)\|_{L^2}^2
\end{equation}
Then $L^2$ convergence of $\tilde u$ implies $\sup_{t\geq 2} \|\tilde u\|_{W^{m, 2}}<\infty$. 

 Suppose $v=\tilde u-\varphi_*$ does not converge to $0$ in $W^{m, 2}$, then there exists $\varepsilon>0,\ t_i\to\infty$ such that $\|v(t_i)\|_{W^{m, 2}}\geq\varepsilon$. 
 By \cite[Proposition B.2]{Bernstein_2017}, the natural embedding $W^{1, 2}(\Sigma, d\tilde \mu)\hookrightarrow L^2(\Sigma, d\tilde \mu)$ is compact. 
 Since $v$ is bounded in $W^{m+1, 2}$, there exists a convergent subsequence of $\{v(t_i)\}$ in $W^{m,2}$. 
 The limit must be $0$, as $v$ converges to $0$ in $L^2$.   This contradicts $\|v(t_i)\|_{W^{m, 2}}\geq\varepsilon$. 
 Therefore, $e^{-\lambda_*t}u$ converges to $\varphi_*$ in $W^{m,2}$ for every $m\in \mathbb N$. On each fixed compact set, the weighted measure $d\tilde\mu$ is comparable to the unweighted one, usual Sobolev embeddings give $C^\infty_{loc}$ convergence.
\end{proof}
\begin{proof}[Proof of Corollary \ref{cor:rigidity}]
   If $u$ is not identically $0$, then Theorem ~\ref{main thm} implies 
\[
\lim_{t\to\infty}e^{kt}\|u\|_{L^2}=\infty,\quad k>-\lambda_*
\]
This contradicts \eqref{cdt2}.
\end{proof}

\section{A local flow with super-exponential decay}\label{sec:local}
In this section, we construct a counterexample to the rigidity under super-exponential decay, under the local graphical assumption. We use a rotationally symmetric shrinker as the base, and find a rotationally symmetric solution to the rescaled equation over an expanding domain. The solution has super-exponential $C^2$ decay and remains nonzero at a point.

\subsection{The base shrinker and the spatial Cauchy problem}

We choose the base shrinker to be the $k=2$ example constructed by Drugan--Kleene \cite[Theorem 4]{drugan2013immersedselfshrinkers}.
This is an immersed shrinker with topology $\R\times\Sn$, rotationally symmetric and asymptotically conical. Let $\gamma(s)=(z(s), r(s))$ be the plane curve that generates the shrinker. The curve has two reflected trumpet ends asymptotic to $r=\pm\sigma z$, for some $\sigma>0$.

Let $\kappa$ be the signed planar curvature, and $\tau, q$ be two supporting functions such that $\gamma=\tau\mathbf T+q\mathbf N$. Then
\[
 \kappa=z'r''-r'z'',\qquad \tau=zz'+rr',\qquad q=-zr'+rz'.
\]
Our construction only needs the following global
bounds, whose proof is in Appendix.

\begin{lemma}\label{lem:geometry}
The profile curve satisfies
\begin{equation}\label{eq:geometry}
 \inf r>0,\qquad
 \sup\bigl(|\kappa|+|\kappa_s|+|q|\bigr)<\infty,\qquad
 |\tau(s)|\leq C(1+|s|).
\end{equation}
\end{lemma}

Solving for the highest spatial derivative, the rescaled mean curvature flow equation is equivalent to
\begin{equation}\label{u eq}
 u_{ss}=\bigl((1-\kappa u)^2+u_s^2\bigr)u_t+G(s,u,u_s),
\end{equation}
where $G(s,0,0)=0$. The formula for $G$ and the derivation of \eqref{u eq} is in Appendix \eqref{app:graph}.

Write $p=u_s$, set
\begin{equation}\label{Y M F def}
 Y=\binom{u}{p},\qquad
 M(s,Y)=\begin{pmatrix}0&0\\(1-\kappa u)^2+p^2&0\end{pmatrix},\qquad
 F(s,Y)=\binom{p}{G(s,u,p)}.
\end{equation}
We are going to solve the Cauchy problem
\begin{equation}\label{Y eq}
 Y_s=M(s,Y)Y_t+F(s,Y),\qquad Y(0,t)=\binom{e^{-t^2}}0.
\end{equation}

For a time series centered at $T$, define the coefficient norm $\|\cdot\|_\rho$ by
\[
 f(t)=\sum_{j\geq0}f_j(t-T)^j,\qquad
 \|f\|_\rho=\sum_{j\geq0}|f_j|\rho^j.
\]
For vectors and matrices sum the component norms. This controls the
maximum on $|t-T|\leq\rho$. We will keep the center implicit when it is fixed, and write $\|f\|_{T,\rho}$ when comparing different centers.
Let $P_N$ be the truncation to degree $N$, $P_Nf=\sum_{i=0}^{N}f_i(t-T)^i$.

\begin{lemma}\label{norm}
The coefficient norm satisfies
\begin{equation}\label{eq:norm-rules}
\begin{gathered}
 \|fg\|_\rho\leq\|f\|_\rho\|g\|_\rho,\qquad
 \|\partial_t f\|_{\rho'}\leq\frac{\|f\|_\rho}{\rho-\rho'},\\
 \|f-P_Nf\|_{\rho'}\leq
 \left(\frac{\rho'}\rho\right)^{N+1}\|f\|_\rho,\qquad 0<\rho'<\rho.
\end{gathered}
\end{equation}
If $[T'-\rho', T'+\rho']\subset[T-\rho, T+\rho]$, then
$\|f\|_{T',\rho'}\leq\|f\|_{T,\rho}$.
\end{lemma}

\begin{lemma}\label{Lip cond}
There are $\delta>0$ and $M_0\geq1$ such that $M,F$ are analytic in $u, p$ on a uniform neighborhood of zero.
For $\|Y\|_\rho,\|Z\|_\rho\leq\delta$,
\begin{equation}\label{eq:Lip}
\begin{aligned}
 \|M(s,Y)\|_\rho&\leq M_0,\\
 \|M(s,Y)-M(s,Z)\|_\rho&\leq M_0\|Y-Z\|_\rho,\\
 \|F(s,Y)\|_\rho&\leq M_0(1+|s|)\|Y\|_\rho,\\
 \|F(s,Y)-F(s,Z)\|_\rho&\leq M_0(1+|s|)\|Y-Z\|_\rho.
\end{aligned}
\end{equation}
\end{lemma}

The proofs of Lemmas~\ref{norm} and \ref{Lip cond} are in the appendix. Fix $\delta,M_0$ from now on and take $a=(256M_0)^{-1}$.

\subsection{Local construction near a large time}

\begin{lemma}\label{local sol}
For every sufficiently large $T$, \eqref{Y eq} has a solution
$Y^T$ on $|s|\leq aT$, $|t-T|<T/16$, smooth in $s$ and analytic in $t$,
such that
\begin{equation}\label{eq:local-bound}
 \sup_{|s|\leq aT}\|Y^T(s,\cdot)\|_{T/16}\leq e^{-T^2/2}.
\end{equation}
\end{lemma}
\begin{proof}
We first give the argument for $s\geq0$, with all norms centered at $T$.

%\paragraph{Uniform estimate on truncated equations.}
Let $Y_N(s, t)=\sum_{j=0}^{N}a_{N,j}(s)(t-T)^j$ be the solution of the truncated equation
\[
 (Y_N)_s=P_N[M(s,Y_N)(Y_N)_t+F(s,Y_N)],\qquad Y_N(0, t)=P_N\binom{e^{-t^2}}{0}
\]
By comparing the coefficients, we get an ODE of the vector $\mathbf{a}_N(s)=(a_{N,j}(s))_{0\leq j\leq N}$. The ODE has a short time solution by the Lipschitz condition of $\mathbf{M, F}$.

Set $\rho(s)=T/8-M_0s$, for almost every $0\leq s\leq aT$ where $\|Y_N\|_{\rho(s)}<\delta$,
\begin{align*}
 \frac d{ds}\|Y_N\|_{\rho(s)}
 &=\sum_{j=0}^\infty\frac d{ds}|a_{N,j}|\rho^j
   +\rho'\sum_{j=1}^\infty j|a_{N,j}|\rho^{j-1}\\
 &\leq\|(Y_N)_s\|_\rho-M_0\|(Y_N)_t\|_\rho\\
 &\leq M_0\|(Y_N)_t\|_\rho+M_0(1+s)\|Y_N\|_\rho-M_0\|(Y_N)_t\|_\rho\\
 &=M_0(1+s)\|Y_N\|_\rho
\end{align*}
By Gronwall inequality, 
\[
\|Y_N\|_{\rho(s)}\leq \|e^{-t^2}\|_{T/8}\exp\left[M_0\left(s+\frac{s^2}{2}\right)\right]
\]
Note that $\|e^{-t^k}\|_\rho=e^{\rho^k}$, $e^{-t^2}=e^{-T^2}e^{-2T(t-T)}e^{-(t-T)^2}$, hence
\[
 \|e^{-t^2}\|_{T/8}\leq e^{-T^2}\|e^{-2T(t-T)}\|_{T/8}\|e^{-(t-T)^2}\|_{T/8}
 = e^{-T^2}e^{2T(T/8)}e^{(T/8)^2}
 =e^{-47T^2/64}
\]
By the choice of $a$, we have $\rho(s)>T/16$ for $0\leq s \leq aT$, hence
\begin{equation}\label{YN estimate}
    \|Y_N\|_{T/16}\leq\|Y_N\|_{\rho(s)}\leq  \exp\left[-\frac{47}{64}T^2+M_0\left(aT+\frac{a^2T^2}{2}\right)\right]\leq e^{-T^2/2},
    \qquad 0\leq s\leq aT
\end{equation}
Choose $T$ large such that $e^{-T^2/2}<\delta/4$. 
The solution $\mathbf{a}_N(s)$ stays in the bounded region $\{ (a_0,\cdots, a_N):\sum_{j=0}^N |a_j|(T/16)^j<\delta/4\}$. Then the ODE continuation theorem can extend the solution to $0\leq s\leq aT$.

%\paragraph{Convergence to the solution.}
The equation of $Y_N$ also gives a bound for $(Y_N)_s$
\[
\|(Y_N)_s\|_{T/32}\leq M_0\|(Y_N)_t\|_{T/32}+M_0(1+s)\|Y_N\|_{T/32}
\leq [\frac{16M_0}{T}+M_0(1+aT)]e^{-T^2/2}
\]
Hence $\{a_{N, j}(s)\}_{N\geq j}$ is uniformly bounded and Lipschitz in $s$. By Ascoli-Arzela theorem and a diagonal argument, $a_{N, j}(s)$ (passing to a subsequence) converges to some $a_j(s)$ uniformly for every $j$.
Define
\[
Y(s, t)=\sum_{j=0}^\infty a_j(s)(t-T)^j.
\]
The power series converges since $\sum_{j=0}^N a_{N,j}(s)(T/16)^j\leq e^{-T^2/2}$ for every $N$. We will show that $Y$ indeed solves \eqref{Y eq}.

For $0<\rho<T/16$ and $J>0$,
\[
\|Y_N-Y\|_\rho
%\leq \sum_{j=0}^{J} |a_{N,j}(s)-a_j(s)|\rho^j+(\frac{16\rho}{T})^{J+1}\|Y_N-Y\|_{T/16}
\leq \sum_{i=0}^{J} |a_{N,j}(s)-a_j(s)|\rho^j+(\frac{16\rho}{T})^{J+1}2e^{-T^2/2}
\]
Let $N\to\infty$, the first term goes to $0$ due to the uniform convergence of $a_{N, j}(s)$. Then let $J\to\infty$, the second term also vanishes, hence $\lim_{N\to\infty} \|Y_N-Y\|_\rho=0$. The convergence is uniform in $s$. Choosing $\eta<T/16$ larger than $\rho$, we have
\[
\|(Y_N-Y)_t\|_{\rho}\leq \frac {\|Y_N-Y\|_{\eta}}{\eta-\rho}\to 0
\]

Set $H_N=M(s,Y_N)(Y_N)_t+F(s,Y_N),\ H=M(s,Y)Y_t+F(s,Y).$
The uniform convergence of $Y_N, (Y_N)_t$ implies $H_N\to H$ uniformly in $s$, in every norm of radius less than $T/16$.
\[
\|H_N-P_N H_N\|_\rho\leq (\frac{\rho}{\eta})^{N+1}\|H_N\|_\eta\to 0,
\quad 0<\rho<\eta<T/16
\]
Therefore, the right hand side of the truncated equation $P_NH_N$ converges to $H$.

Passing limit to the integrated equation
\[
 Y_N(s)=P_N\binom{e^{-t^2}}{0}+\int_0^sP_NH_N(\sigma)\,d\sigma
\]
we get
\[
Y(s)=\binom{e^{-t^2}}{0}+\int_0^sH(\sigma)\,d\sigma
\]
Hence $Y$ solves our Cauchy problem \eqref{Y eq}, and satisfies the estimate \eqref{eq:local-bound}.

For $s<0$, set $\rho(s)$ by $\rho(s)=T/8+M_0s$, the above argument still works. The equation determines matching derivatives from the same initial data at $s=0$, so the two parts join smoothly.
\end{proof}

\subsection{Agreement on overlaps and solution for all time}

The local solutions have the same initial data at $s=0$. The next lemma shows that they agree away from that line on smaller time intervals. 

\begin{lemma}\label{overlap}
For sufficiently large $T,T'$, the solutions $Y^T,Y^{T'}$ agree
wherever the rectangles
\[
 \{|s|<aT,\ |t-T|<T/64\},\qquad
 \{|s|<aT',\ |t-T'|<T'/64\}
\]
overlap.
\end{lemma}
\begin{proof}
Let $W=Y^T-Y^{T'}$, then $W$ satisfies the equation
\[
W_s=\mathbf{M}(Y^T,s)W_t+(\mathbf{M}(Y^{T},s)-\mathbf{M}(Y^{T'},s))Y^{T'}_t+\mathbf{F}(Y^T,s)-\mathbf{F}(Y^{T'}, s)
%,\qquad W(0, t)=\mathbf{0}
\]
Without loss of generality, assume $T'>T$. Take $t_*\in[T'-T'/64, T+T/64]$, then $ $ $[t_*-T/32, t_*+T/32]\subset[T'-T'/16, T+T/16]$. 
We showed that $\|Y^T\|_{T, T/16}, \|Y^{T'}\|_{T',T'/16}\leq \delta/4$. By Lemma ~\ref{norm}, $\|Y^T\|_{t_*, T/32}, \|Y^{T'}\|_{t_*,T/32}\leq \delta/4$, so $\|W\|_{t_*,T/32}\leq \delta/2$.
\[
\|Y^{T'}_t\|_{t_*,T/64}\leq \frac{64\|Y^{T'}\|_{t_*, T/32}}{T}\leq \frac{16\delta}{T}
\]

The following norms are all centered at $t_*$.
Let $\rho(s)=T/64-M_0s$. For $0\leq s<aT$, $\rho(s)>T/128$. Similar to the previous computation, we have
\begin{align*}
 \frac d{ds}\|W\|_{\rho(s)} 
 &\leq\|W_s\|_\rho-\rho'\|W_t\|_\rho\\
 &\leq M_0\|W_t\|_\rho+M_0\frac{16\delta}{T}\|W\|_\rho+M_0(1+|s|)\|W\|_\rho-M_0\|W_t\|_\rho\\
 &=M_0(\frac{16\delta}{T}+1+|s|)\|W\|_\rho
\end{align*}
Since $W(0, t)=\mathbf{0}$, Gronwall inequality implies $\|W\|_{\rho(s)}\equiv\mathbf{0}$. In particular, $W(s, t_*)=0$. Replacing $s$ by $-s$, above argument also works for negative parameter.
The choice of $t_*$ is arbitrary, hence $Y^{[T]}=Y^{[T']}$ on the overlap.
\end{proof}

\begin{proof}[Proof of Theorem~\ref{counterexample}]
Choose large time $T_0$ such that previous lemmas hold for $T\geq T_0$. We can construct solution $Y$ for all time beyond $T_0$.

Let $T_n=T_0+n$, define $Y(s, t)=Y^{[T_n]} (s,t) $ in $\Omega_n=\{(s,t):|s|<aT_n, |t-T_n|< T_n/64\}$. For each $(s, t)$ such that $t>T_0, |s|<at$, there are only finitely many $T_n$ satisfy $(s, t)\in\Omega_n $, and Lemma ~\ref{overlap} ensures $Y^{[T_n]}(s, t)$ takes the same value for each $n$. Hence $Y$ is well-defined in $\{(s, t): |s|<at, t>T_0\}\subset\cup_n \Omega_n$. $Y$ is nonzero since $Y(0, t)=(e^{-t^2}, 0)$.

Next, we check the super-exponential decay condition\eqref{exp decay}.
For $t>T_0, |s|<at$, consider the norm centered at $t$.
\[
\|Y_t\|_{t/32}\leq \frac{32\|Y\|_{t/16}}{t}\leq \frac{32}{T_0}e^{-t^2/2}
\]
\[
\|Y_s\|_{t/32}\leq M_0\|Y_t\|_{t/32}+M_0(1+|s|)\|Y\|_{t/32}\leq M_0(\frac{32}{T_0}e^{-t^2/2}+1+at)
\]
Recall $Y=(u, u_s)$, the first two spatial derivatives of $u$ is controlled by $Y, Y_s$
\[
\sup_{|s|<at}|u|+|u_s|+|u_{ss}|\leq \|Y\|_{t/32}+\|Y_s\|_{t/32}\leq M_0\max\{a, 1+32/T_0\}(1+t)e^{-t^2/2}
\]
Under the metric $g_\Sigma=ds^2+r^2g_{\Sn}$, the eigenvalues of $\nabla^2_\Sigma u$ are $u_{ss}$ and
$(r'/r)u_s$ with multiplicity $n-1$. The covariant derivatives of $u$ are
\[
|\nabla_\Sigma u|=|u_s|,\quad |\nabla_\Sigma^2 u|^2=|u_{ss}|^2+(n-1)(\frac{r'}{r})^2|u_s|^2
\]
Since $r'/r$ is bounded, there exists $C>0$ such that
\[
\sup_{|s|<at}|u|+|\nabla_\Sigma u|+|\nabla^2_\Sigma u|\leq Ce^{-t^2/4}
\]
For any compact set $V\subset \Sigma$, there exists $T_0<T_V<\infty$ such that $V \subset [-aT_V, aT_V]\times \Sn$. 
\[
\|u(\cdot,t)\|_{C^2(V)}\leq \sup_{|s|<aT_V}(|u|+|\nabla_\Sigma u|+|\nabla^2_\Sigma u|)(s,t)\leq Ce^{-t^2/4},\quad t\geq T_V.
\]
Thus, $\lim_{t\to \infty} e^{kt}\|u(\cdot , t)\|_{C^2(V)}=0$. 
\end{proof}

\appendix
\section{Calculations of the profile curve and coefficient bounds}\label{sec:appendix}

We will verify Lemmas~\ref{lem:geometry}, \ref{norm}, and \ref{Lip cond}, and derive the exact normal-graph equation. 

\begin{lemma}[Lemma~\ref{lem:geometry}]
    The profile curve satisfies
\begin{equation*}
 \inf r>0,\qquad
 \sup\bigl(|\kappa|+|\kappa_s|+|q|\bigr)<\infty,\qquad
 |\tau(s)|\leq C(1+|s|).
\end{equation*}
\end{lemma}
\begin{proof}
Let $\mathbf T, \mathbf N$ be the planar tangent and normal vector, then
\[
 \mathbf T=(z',r'),\qquad \mathbf N=(-r',z'),\qquad
 \mathbf T'=\kappa\mathbf N,\qquad z'^2+r'^2=1.
\]
The hypersurface normal is $\nu=(-r',z'\omega)$.
The shrinker equation and Frenet formulas give
\begin{equation}\label{eq:profile}
 -\kappa+(n-1)\frac{z'}r=\frac q2,\qquad
 \tau'=1+\kappa q,\qquad q'=-\kappa\tau.
\end{equation}
On each end $r$ is a function $f$ of $z$, the trumpet estimates
\cite[Theorem 3]{Kleene_2014} imply there exists $\sigma>0$, such that
\[
 f(z)=\sigma|z|+O(|z|^{-1}),\qquad
 f_z(z)=\sigma\operatorname{sgn}(z)+O(|z|^{-2})
\]
Therefore, the arclength is comparable to $|z|$.
Up to orientation, $q=(f-zf_z)/\sqrt{1+f_z^2}=O(|s|^{-1})$,
hence $\kappa=O(|s|^{-1})$.
Differentiating the first identity in \eqref{eq:profile} gives
\[
 \kappa_s=\left(\frac\tau2-(n-1)\frac{r'}r\right)\kappa
              -(n-1)\frac{z'r'}{r^2}.
\]
Unit speed implies $|\tau(s)|\leq |\gamma(s)|\leq|\gamma(0)|+|s|$, so $\kappa_s$ is
bounded on both ends. The remaining compact part lies in
$r>0$. This proves all the bounds.
\end{proof}

\begin{figure}[!htb]
 \centering
 \includegraphics[width=.42\textwidth]{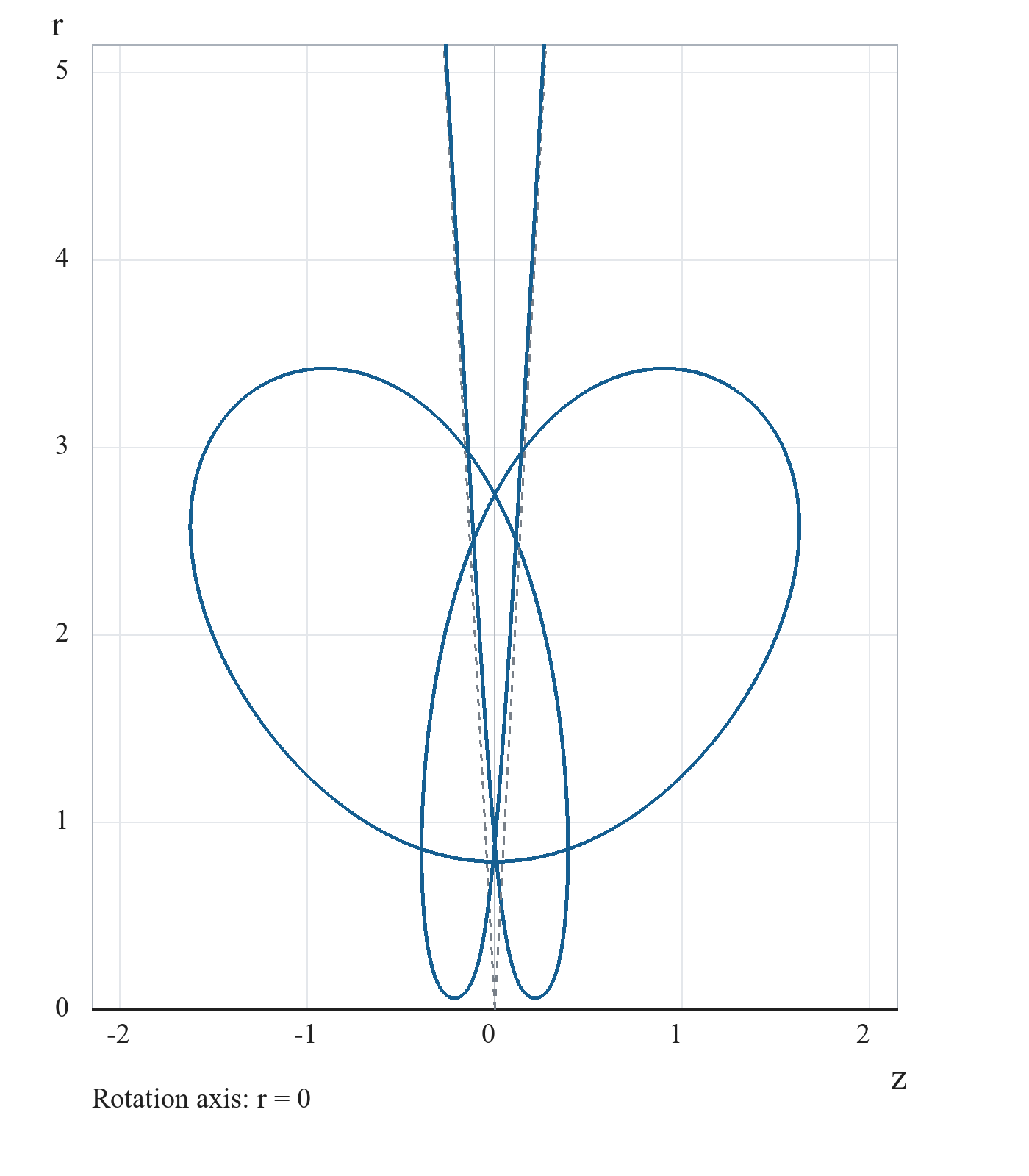}
 \caption{Numerical illustration of the two-ended $k=2$ profile curve. Dashed lines indicate the asymptotic cone rays.}
 \label{fig:profile}
\end{figure}

\subsection{The exact normal-graph equation}\label{app:graph}
Define the following geometric factors:
\[
 p=u_s,\qquad h=1-\kappa u,\qquad w^2=h^2+p^2,\qquad
 \ell=r+uz'
\]
The perturbed planar curve $\gamma+u\mathbf N$ has unit tangent
$\mathbf T_u=(h\mathbf T+p\mathbf N)/w$ and unit normal $\mathbf N_u=(-p\mathbf T+h\mathbf N)/w$.
Frenet formula gives the signed curvature
\[
 \kappa_u=\frac{\ip{\mathbf T_u'}{\mathbf N_u}}{w}=\frac{\kappa h^2+h u_{ss}+\kappa_s up+2\kappa p^2}{w^3}
\]
The principal curvatures of the perturbed surface are $-\kappa_u$, and
$(hz'-pr')/(w\ell)$ with multiplicity $n-1$. Hence
\[
 H_u=-\kappa_u+(n-1)\frac{hz'-pr'}{w\ell},\qquad
 \ip{X_u}{\nu_u}=\frac{h(q+u)-p\tau}{w},\qquad
 \ip{\partial_tX_u}{\nu_u}=\frac hw u_t.
\]
Substituting this into the rescaled equation \eqref{eq:rmcf}, we get
\begin{equation}\label{eq:graph}
 u_t=\frac{u_{ss}}{w^2}
 +\frac{\kappa h^2+2\kappa p^2+\kappa_s up}{hw^2}+(n-1)\left(\frac{pr'}{h\ell}-\frac{z'}\ell\right)
 -\frac{p\tau}{2h}+\frac{q+u}{2}.
\end{equation}
Solving for $u_{ss}$ and canceling the constant terms by
\eqref{eq:profile} gives
\[
 u_{ss}=w^2u_t+G(s,u,p)
\]
where
\begin{equation}\label{eq:G}
\begin{split}
 G(s,u,p)={}&-\kappa(h-1)-\frac{2\kappa p^2+\kappa_s up}{h}
 -(n-1)\frac{w^2pr'}{h\ell}\\
 &+(n-1)z'\left(\frac{w^2}{\ell}-\frac1r\right)
 -\frac q2(w^2-1)-\frac{w^2u}{2}
 +\frac{\tau w^2p}{2h}.
\end{split}
\end{equation}
Every term vanishes when $(u,p)=(0,0)$. This proves \eqref{u eq}.

\subsection{The coefficient norm and Lipschitz conditions}
\begin{lemma}[Lemma~\ref{norm}]
The coefficient norm satisfies:
    \begin{equation*}
\begin{gathered}
 \|fg\|_\rho\leq\|f\|_\rho\|g\|_\rho,\qquad
 \|\partial_t f\|_{\rho'}\leq\frac{\|f\|_\rho}{\rho-\rho'},\\
 \|f-P_Nf\|_{\rho'}\leq
 \left(\frac{\rho'}\rho\right)^{N+1}\|f\|_\rho,\qquad 0<\rho'<\rho.
\end{gathered}
\end{equation*}
If $[T'-\rho', T'+\rho']\subset[T-\rho, T+\rho]$, then
$\|f\|_{T',\rho'}\leq\|f\|_{T,\rho}$.
\end{lemma}
\begin{proof}
Multiplication of series gives the product bound. 

For $0<\rho'<\rho, \ j\geq 1$, we have $ \rho^j>\rho^j-(\rho')^{j}\geq(\rho-\rho')j(\rho')^{j-1}$.

\[
\|\partial_t f\|_{\rho'}=
\sum_{j\geq1} j(\rho')^{j-1}|f_j|\leq
\sum_{j\geq1}\frac{\rho^j|f_j|}{\rho-\rho'}=\frac{\| f\|_\rho}{\rho-\rho'}
\]
\[
\|f-f_N\|_{\rho'}=\sum_{j\geq N+1}|f_j|\rho^j(\frac{\rho'}{\rho})^{j}\leq \sum_{j\geq N+1}|f_j|\rho^j(\frac{\rho'}{\rho})^{N+1}
= (\frac{\rho'}{\rho})^{N+1}\|f\|_\rho
\]
Finally expand $(t-T)^j=((t-T')+(T'-T))^j$ by the binomial formula.
\[
f(t)=\sum_{j\geq 0}f_j\sum_{k=0}^j\binom{j}{k}(t-T')^k(T'-T)^{j-k}=\sum_{k\geq 0}\big(\sum_{j\geq k}f_j\binom{j}{k}(T'-T)^{j-k}\big)(t-T')^k
\]
\[
 \|f\|_{T',\rho'}\leq\sum_{k\geq 0}\big(\sum_{j\geq k}|f_j|\binom{j}{k}|T'-T|^{j-k}\big)(\rho')^k
 =\sum_{j\geq0}|f_j|(|T'-T|+\rho')^j\leq\|f\|_{T,\rho}.
\]
\end{proof}

\begin{lemma}[Lemma~\ref{Lip cond}]
There are $\delta>0$ and $M_0\geq1$ such that $M,F$ are analytic in $u, p$ on a uniform neighborhood of zero.
For $\|Y\|_\rho,\|Z\|_\rho\leq\delta$,
\begin{equation*}
\begin{aligned}
 \|M(s,Y)\|_\rho&\leq M_0,\\
 \|M(s,Y)-M(s,Z)\|_\rho&\leq M_0\|Y-Z\|_\rho,\\
 \|F(s,Y)\|_\rho&\leq M_0(1+|s|)\|Y\|_\rho,\\
 \|F(s,Y)-F(s,Z)\|_\rho&\leq M_0(1+|s|)\|Y-Z\|_\rho.
\end{aligned}
\end{equation*}
\end{lemma}
\begin{proof}
Choose
\[
 \delta=\frac1{4(1+\sup|\kappa|+\sup r^{-1})}.
\]
For $\|Y\|_\rho=\|u\|_\rho+\|p\|_\rho\leq\delta$, geometric series give
\begin{align*}
    &\|h^{-1}\|_\rho=\|\sum_{j=0}^\infty (\kappa u)^j\|_\rho\leq \sum_{j=0}^\infty (\sup |\kappa|\delta) ^j=\frac{4}{3},\qquad
\|\ell^{-1}\|_\rho\leq \frac{1}{r}\sum_{j=0}^\infty (\sup |r^{-1}|\delta) ^j\leq\frac4{3r},\\
 &\|w^2\|_\rho\leq(5/4)^2+\delta^2<2.
\end{align*}
These series also show analyticity in $u, p$ on a uniform small ball centered at $0$. In \eqref{eq:G}, use
\[
 w^2-1=-2\kappa u+\kappa^2u^2+p^2,\qquad
 \frac{w^2}{\ell}-\frac1r
 =\frac{w^2-1}{\ell}-\frac{uz'}{r\ell}.
\]
Then every single term of $G$ contains a graph factor $u$ or $p$.
Lemma~\ref{lem:geometry} bounds all geometric coefficients except
$\tau$, which grows at most linearly on $|s|$. Therefore
\[
 \|M(s,Y)\|_\rho\leq C,\qquad
 \|F(s,Y)\|_\rho\leq C(1+|s|)\|Y\|_\rho.
\]
Subtract products one factor at a time to bound the differences.
For inverse factors,
\[
\| h_Y^{-1}-h_Z^{-1}\|_\rho=\|h_Y^{-1}(h_Z-h_Y)h_Z^{-1}\|_\rho
\leq\frac{16}{9}\|\kappa(u_Y-u_Z)\|_\rho\leq \frac{4}{9\delta}\|Y-Z\|_\rho
\]
and similarly for $\ell^{-1}$. The same bounds give
\[
 \|M(s,Y)-M(s,Z)\|_\rho\leq C\|Y-Z\|_\rho,\qquad
 \|F(s,Y)-F(s,Z)\|_\rho\leq C(1+|s|)\|Y-Z\|_\rho.
\]
Choose $M_0\geq1$ larger than the above constants, then all the estimates are satisfied.
\end{proof}

\newpage
\bibliographystyle{plain}
\bibliography{references}
\end{document}